\documentclass{amsproc}
\usepackage{amssymb, amsmath}
\usepackage{xypic}
\xyoption{all}

\swapnumbers
\theoremstyle{theorem}
\newtheorem{theorem}{Theorem}[section]
\newtheorem{lemma}[theorem]{Lemma}
\newtheorem{proposition}[theorem]{Proposition}
\newtheorem{corollary}[theorem]{Corollary}

\theoremstyle{definition}
\newtheorem{definition}[theorem]{Definition}
\newtheorem{example}[theorem]{Example}
\newtheorem{question}[theorem]{Question}

\newtheorem{point}[theorem]{}

\theoremstyle{remark}

\def\C{{\mathcal {C}}}
\def\D{{\mathcal {D}}}

\def\Z{{\mathbf {Z}}}

\def\Arrows{\it {Arrows}}

\def\rmono{\rto|<\hole|<<\ahook}

\def\umono{\ar@{_{(}->}[u]}
\def\uumono{\ar@{_{(}->}[uu]}

\def\lmono{\ar@{_{(}->}[l]}
\def\llmono{\ar@{_{(}->}[ll]}

\def\rrepi{\rrto|>>\tip}

\def\ra{\rightarrow}
\def\la{\leftarrow}
\def\xra{\xrightarrow}
\def\xla{\xleftarrow}

\def\mono{\hookrightarrow}

\begin{document}

\title{Homotopy pull-back squares up to localization}

\author{Wojciech Chach\'olski}
\address{K.T.H.\\ Matematik\\ S-10044 Stockholm,
Sweden}
\email{wojtek@math.kth.se}

\author{Wolfgang Pitsch}
\address{Universitat Aut\`onoma de Barcelona \\ Departament de Matem\`atiques\\
E-08193 Bellaterra, Spain}
\email{pitsch@mat.uab.es}

\author{J\'er\^ome Scherer}
\address{Universitat Aut\`onoma de Barcelona \\ Departament de Matem\`atiques\
\
E-08193 Bellaterra, Spain}
\email{jscherer@mat.uab.es}

\thanks{The first author is supported in part by  Vetenskapsr\aa
det grant 2001-4296 and G\"oran Gustafssons Stiftelse, the second and third
authors by MEC grant MTM2004-06686 and  by the program Ram\'on y Cajal, MEC, Spa
in}
\date{December 31, 2004.}


\begin{abstract}
We characterize the class of homotopy pull-back squares by means of elementary
closure properties.
The so called Puppe theorem which identifies the homotopy fiber of certain maps 
constructed as homotopy colimits
is a straightforward consequence. Likewise we characterize the class of squares 
which are homotopy pull-backs ``up to Bousfield localization". This yields
a generalization of Puppe's theorem which allows to identify the homotopy type
of the localized homotopy fiber. When the localization functor is homological
localization this is one of the key
ingredients in the group completion theorem.
\end{abstract}

\maketitle

\section{Introduction}
In topology it is convenient to think about a continuous family of spaces as  a 
map whose fibers constitute the family. The homotopy fiber of this map is then an
important invariant of the family, but in general it is difficult to say anything about
this invariant. However if the ``transition functions" of the map  preserve some property, it is
often the case that the same property is inherited by the homotopy fiber.
The classical example is given by the so called Puppe theorem \cite{MR51:1808} that says that if all the members of
the family have the same homotopy type (the transition functions are weak equivalences), then the
homotopy fiber has this homotopy type too. This is also  the central idea in (generalized) Quillen's group
completion theorem, as exposed by McDuff-Segal~\cite{MR53:6547}, Jardine~\cite{MR90j:55029},
Tillman~\cite{MR99k:57036}, see also Adams' book \cite{MR80d:55001}. In their setting the members of
the family have the same integral homology type (the transition functions are $\text{H}\mathbf{Z}$-isomorphisms).
The statement asserts  then that the homotopy fiber shares the same integral homology type as well. This was
used to compute the homology of the group completion of certain topological monoids: a celebrated
consequence is the Baratt--Priddy theorem \cite{MR47:3489} which identifies
$B\Sigma_\infty^+ \times \Z$ with $QS^0$.

The aim of this paper is to generalize these results to the case when the members of a continuous family of spaces
have the same homotopy type after Bousfield localization with respect to a given map. It has been surprising to us that this
statement follows from properly reformulating the classical Puppe theorem together with general properties of localizations
of spaces.

To state our  theorem it turns out that it is more appropriate to work in the category $\Arrows$ of maps of
spaces (see Section~\ref{pt arrows}) rather than in $Spaces$. Studying homotopy 
fibers of continuous families
can then be translated into investigating  properties of certain classes of morphisms in $\Arrows$:

\begin{definition}\label{def dist}
A class $\C$ of morphisms  in  \mbox{$\Arrows$} is called
\emph{distinguished} if:
\begin{enumerate}
\item  Weak equivalences belong to $\C$.
\item Let $\phi:f\ra g$ and $\psi:g\ra h$ be morphisms. Assume that  either $\psi$
or $\phi$ is a weak equivalence. Then if two out of $\phi$, $\psi$,  $\psi\phi$ 
belong to $\C$, then so does the third.
\item  Let $\phi:f\ra g$ and $\psi:g\ra h$ be morphisms. If $\psi$ and $\psi\phi$ belong to $\C$, then so does $\phi$.
\item If  $F:I\ra \mbox{$\Arrows$}$   sends
any morphism in $I$ to  a morphism in $\C$, then, for every
$i\in I$,  $F(i)\ra \text{hocolim}_{I}F$ belongs to $\C$.
\end{enumerate} 
\end{definition}

At first it might seem pointless to consider such collections since for example 
the category $Spaces$
has only one distinguished class that consists of all spaces. The key observation is that there are much more
interesting distinguished collections in  $\Arrows$. For example our  Puppe theorem can  now be formulated as follows, it
identifies the class of homotopy pull-backs. 
\medskip

\ref{thm mainpullback} T\begin{footnotesize}HEOREM\end{footnotesize}.
{\it The collection of homotopy pull-backs is the smallest distinguished class of morphisms in \mbox{$\Arrows$}.}

\medskip

As in $Spaces$, Bousfield localization also exist in $\Arrows$.  Although it is 
very hard to identify these localization
functors with respect to an arbitrary morphism, an explicit description can be given for  $L_{\phi}$ when
$\phi=(u,\text{id}_{\Delta[0]})$ is a morphism between two maps collapsing a space to a point (see Section~\ref{pt arrows}).
In this case $L_{\phi}$ is written $L_{u}$ and coincides with the ``fiberwise" application of the Bousfield localization
$L_{u}$ of spaces (see Section~\ref{sec fiblocal}). We say that a morphism
$\psi:f\ra g$ in $\Arrows$ is an $L_{u}$-homotopy pull-back if $\psi$ induces
$L_{u}$-equivalences of homotopy fibers of $f$ and $g$ (see Definition~\ref{def h-homotopypb}) or equivalently if
$L_{u}\psi: L_{u} f \rightarrow L_{u} g$ is a homotopy pull-back (see Proposition~\ref{prop amenable}).
Our main theorem can be now stated as follows:

\medskip

\ref{thm luhpullback} T\begin{footnotesize}HEOREM\end{footnotesize}.
{\it The collection of $L_{u}$-homotopy pull-backs is the smallest distinguished class
containing all $L_{u}$-equivalences.}

\medskip

For example let us choose  a map $u$ for which $L_{u}$ coincides with the localization with respect to a chosen
homology theory. Let $F,G:I\ra Spaces$  be functors and $\pi:F\ra G$ be a natural transformation. Assume that  for any
morphism $\alpha\in I$,  the commutative square:
\[\xymatrix{
F(i)\rto^{F(\alpha)}\dto_{\pi_{i}} & F(j)\dto^{\pi_{j}}\\
G(i)\rto^{G(\alpha)} & G(j)
}\]
induces a homology isomorphism between the homotopy fibers of $\pi_{i}$ and $\pi_{j}$
(i.e. this square is  an $L_{u}$-homotopy pull-back). Then, since $L_{u}$-homotopy pull-backs form a
distinguished collection, according to condition (4) of Definition~\ref{def dist}, the homotopy fibers
of $\text{hocolim}_{I}\pi$ have the same homology type as the homotopy fibers of
 $\pi_{i}$ for an appropriate $i$.

\medskip

{\bf Acknowledgments}: It was Emmanuel Dror Farjoun who suggested to view group completion as a
fiberwise localization statement when we explained him the result in \cite{MR2061572}.
Unlike some other participants we were not able to prove our main theorem while hiking
around Arolla. We thank the G\"oran Gustafssons Stiftelsethe for support and KTH in Stockholm for its hospitality.

\section{The category of maps of spaces}
In this section we deal with combinatorics and geometry of simplicial sets.
We focus particularly on geometrical  properties of push-outs and pull-backs.
\begin{point}\label{pt arrows}
The symbol \mbox{$\Arrows$} denotes the category whose objects are
maps in $Spaces$ and morphisms are commutative squares. Explicitly
a morphism $\phi:f\ra g$ in \mbox{$\Arrows$} is given by a pair
 $\phi=(\phi_{0},\phi_{1})$ of maps for which the following square commutes:
\[\xymatrix{
X\rto^{\phi_{0}}\dto_{f} & A\dto^{g}\\
Y\rto^{\phi_{1}} & B
}\]
\end{point}
\begin{point}
There are two forgetful functors $D,R:\mbox{$\Arrows$}\ra Spaces$
which assign to a map $f:X\ra Y$ its domain $Df:=X$ and its range
$Rf:=Y$. So for the morphism $\phi$,  $D\phi=\phi_{0}$ and
$R\phi=\phi_{1}$, which should not be mistaken for the domain $f$ and the range 
$g$ of  $\phi$ considered as a morphism in $\Arrows$.
A functor $F:I\ra \mbox{$\Arrows$}$ can be
identified with a natural transformation $DF\ra RF$ between
functors with values in  $Spaces$, denoted by $\pi F$. By the universal
properties  $\text{colim}_{I}F$ and $\text{lim}_{I}F$ are naturally isomorphic respectively to the maps
$\text{colim}_{I}\pi F:\text{colim}_{I}DF\ra \text{colim}_{I}RF$ and
$\text{lim}_{I}\pi F:\text{lim}_{I}DF\ra \text{lim}_{I}RF$.
\end{point}

\begin{point}\label{pt prop}
A morphism $\phi:f\ra g$  in \mbox{$\Arrows$} is called a
pull-back or a push-out if the corresponding square is so in
$Spaces$. It is called a monomorphism if both $D\phi$ and $R\phi$
are so in $Spaces$.
\end{point}

Here is a list of basic properties of pull-backs and push-outs of
spaces. One way of proving them is to show that they are true for
sets and are preserved by functor categories, thus they remain
valid for $Spaces$, \mbox{$\Arrows$}, etc. The first property is
the classical ``two out of three" property, similar statements can
be found for example in \cite[Proposition~1.8]{MR93i:55015}.

\begin{lemma}\label{lemma 2outof3}
Let $\phi:f\ra g$ and $\psi:g\ra h$ be two morphisms in
\mbox{$\Arrows$}.

\noindent \emph{(1)} Assume that $\psi$ (respectively $\phi$) is a
pull-back (respectively a push-out). Then  $\phi$
 (respectively $\psi$) is a pull-back (respectively a push-out) if and only if
$\psi\phi$ is so.

\noindent \emph{(2)} Assume that $\phi$ is a pull-back and
$R\phi:Rf\ra Rg$ is an epimorphism. Then $\psi$ is a pull-back if
and only if $\psi\phi$ is so.

\noindent \emph{(3)} Assume that $\psi$ is a push-out and $f$,
$g$, $h$ are monomorphisms. Then $\psi\phi$ is a push-out if and
only if $\phi$ is so.\hfill{\qed}
\end{lemma}

The following examples show that the additional assumptions in points (2) and (3) above are essential.

\begin{example}\label{ex not23}
Let $S^0$ denote the boundary of $\Delta[1]$. Here is a diagram in which the left square and the outer square are pull-backs,
but the right square is not:
\[\xymatrix{\emptyset \dto \rto & \Delta[0] \dto^{d^1} \rto &\Delta[0] \dto^{d^1} \\
\Delta[0] \rto^{d^0} & S^0 \rto^{d^0 s^0} & S^0
}
\]

Here is a diagram in which the right square and the outer square are push-outs, 
but the left square is not:
\[\xymatrix{  \Delta[0]  \dto \rto & S^0 \dto \rto &  \Delta[0] \dto \\
 \Delta[0] \rto &   \Delta[0]  \rto &   \Delta[0]
}
\]
\end{example}

The next three properties are occurrences of push-out squares
being at the same time pull-back squares. Consider a commutative
diagram in \mbox{$\Arrows$}:
\begin{equation}\label{diag pushpull}
\diagram
f\rto^{\phi}\dto_{\pi} &g\dto^{\mu}\\
h\rto^{\psi} & k
\enddiagram
\end{equation}

\begin{lemma}\label{lemma pushoutispullback}
\hfill{}

\noindent \emph{(1)} If $\phi:f\ra g$ is a monomorphism and a
push-out, then it is a pull-back.

\noindent \emph{(2)} Assume that the square~{\rm(\ref{diag pushpull})}
is a push-out. If $\phi$ is a monomorphism, then so is $\psi$ and
this square is also a pull-back.

\noindent \emph{(3)} Assume that the square~{\rm(\ref{diag pushpull})}
is a pull-back. If $\psi$ is a pull-back (respectively a
monomorphism), then so is $\phi$. In particular if $\psi$ is a
push-out and a monomorphism, then $\phi$ is a monomorphism and a
pull-back. \hfill{\qed}
\end{lemma}

Finally we state two different ``cube theorems", named in analogy
with Mather's theorem \cite{MR595633}. They say that sometimes
push-outs do commute with pull-backs.

\begin{lemma}\label{lemma cube}
Assume that in the square~{\rm(\ref{diag pushpull})} the morphisms $\phi$ and $\psi$ (respectively $\phi$,
$\psi$, $\pi$, and $\mu$) are pull-backs. Then {\rm(\ref{diag pushpull})} is a pull-back (respectively a
push-out) square if and only if the range square:
\[\xymatrix{
Rf\rto^{R\phi} \dto_{R\pi}& Rg\dto^{R\mu}\\
Rh\rto^{R\psi} & Rk }\] is a pull-back (respectively a push-out)
of spaces.  \hfill{\qed}
\end{lemma}

This statement, which will be referred to as the cube lemma, has
the following extension to \mbox{$\Arrows$}. We call it the
hypercube lemma. Consider a commutative cube in \mbox{$\Arrows$}:
\[\xymatrix{
 &\bar{f}\rrto\dlto\ddto|\hole& & \bar{g}\dlto\ddto\\
 \bar{h}\rrto\ddto & & \bar{k}\ddto\\
 & f\rrto|\hole\dlto & & g\dlto\\
 h\rrto & & k
}\]

\begin{lemma}\label{lemma hypercube}
In the above cube assume that the squares $(\bar{f},\bar{h},f,h)$,
$(\bar{f},\bar{g},f,g)$, $(\bar {h},\bar{k},h,k)$, and
$(\bar{g},\bar{k},g,k)$ are pull-backs and that the square
$(f,h,g,k)$ is a push-out. Then the square
$(\bar{f},\bar{h},\bar{g},\bar{k})$ is also a push-out.   \hfill{\qed}
\end{lemma}

\begin{point}
Fix a morphism $\phi:f\ra g$ in \mbox{$\Arrows$}. There are two
natural operations one can perform. First, given a morphism
$\tau \rightarrow g$ one can \emph{pull it back} along $\phi$ and define
$\phi^{\ast}\tau\ra f$ to be the morphism that fits into the
following pull-back square in \mbox{$\Arrows$}:
\[\xymatrix{\phi^{\ast}\tau\rto\dto & \tau\dto\\
f\rto^{\phi}& g }\]
In general $\phi^{\ast}\tau\ra f$ is not a pull-back, but
it is so whenever $\tau\ra g$ is a pull-back.

Second for any $\sigma\ra f$, define the \emph{push-forward}
$\phi_{\ast}\sigma\ra g$ to be a pull-back that fits into a
commutative square in \mbox{$\Arrows$} of the form:
\[
\xymatrix{
\sigma\rto\dto & \phi_{\ast}\sigma\dto\\
f\rto^{\phi} & g }
\]
and is initial with respect to this property. Explicitly
$\phi_{\ast}\sigma\ra g$ is given by the following pull-back
square in $Spaces$:
\[
\xymatrix{
A\rrto \dto_{\phi_{\ast}\sigma} & & Dg\dto^{g}\\
R\sigma\rto & Rf\rto^{R\phi} & Rg }
\]
This construction shows that the push-forward is functorial. Note
that $\phi$ is a pull-back if and only if, for any pull-back
$\sigma\ra f$, the morphism $\sigma\ra \phi_{\ast}\sigma$ is an
isomorphism. By definition the push-forward
$\phi_{\ast}\sigma\ra g$ is always a pull-back.

Observe that there is a natural morphism
$\sigma\ra \phi^{\ast}\phi_{\ast}\sigma $, which is a pull-back if
$\sigma\ra f$ is so. For any pull-back $\tau\ra g$, there is also a natural
pull-back morphism  $\phi_{\ast}\phi^{\ast}\tau\ra \tau$.
\end{point}

\begin{proposition}\label{prop commut}
Assume that the range square of the square~{\rm(\ref{diag pushpull})}
is a pull-back of spaces. For any pull-back $\sigma\ra h$, the induced morphism
$\phi_{\ast}\pi^{\ast}\sigma\ra \mu^{\ast}\psi_{\ast}\sigma$ is then an
isomorphism.
\end{proposition}

\begin{proof}
To prove the proposition we need to show that the following square
is a pull-back in \mbox{$\Arrows$}:
\[\xymatrix{
\phi_{\ast}\pi^{\ast}\sigma\rto\dto & g\dto^{\mu}\\
\psi_{\ast}\sigma\rto & k }\]
Since the horizontal morphisms in
this square are pull-backs, according to the cube
Lemma~\ref{lemma cube}, it is enough to show that this square is a pull-back on the
range level. On the range level, this is the outer square of:
\[\xymatrix{
R(\pi^{\ast}\sigma)\rto\dto & Rf\rto^{R\phi}\dto^{R\pi} &
Rg\dto^{R\mu}
\\
R\sigma\rto & Rh\rto^{R\psi} & Rk }\]
As the left and right squares of this diagram are pull-backs, then so is the outer one,
proving the proposition.
\end{proof}

\begin{proposition}\label{prop pushfiber}
Assume that the square~{\rm(\ref{diag pushpull})} is a push-out
and either $\phi$ is a monomorphism or $\pi$ is a monomorphism and
a pull-back. Then for any pull-back $\sigma \ra h$ the following
is a push-out square:
 \[\xymatrix{
\pi^{\ast}\sigma\rto\dto & \phi_{\ast}\pi^{\ast}\sigma\dto\\
\sigma\rto & \psi_{\ast}\sigma
}\]
\end{proposition}

\begin{proof}
By Lemma~\ref{lemma pushoutispullback}.(1) the
square~{\rm(\ref{diag pushpull})} is also a pull-back. Hence
according to Proposition~\ref{prop commut}, the morphism
$\phi_{\ast}\pi^{\ast}\sigma\ra \mu^{\ast}\psi_{\ast }\sigma$ is
an isomorphism. Consider next the following commutative diagram in
\mbox{$\Arrows$}:
\begin{equation}\label{biggd} 
\diagram
 & \pi^{\ast}\sigma\rto\dlto\ddrto|(0.34)\hole|\hole &  \phi^{\ast}\phi_{\ast}
 \pi^{\ast}\sigma\rrto\dlto\ddto|\hole & & \mu^{\ast}\psi_{\ast}\sigma\dlto\ddto
 &
 \hspace{-10mm} =\phi_{\ast}\pi^{\ast}\sigma\\
\sigma\rto\ddrto &\psi^{\ast}\psi_{\ast}\sigma\rrto \ddto& &\psi_{\ast}\sigma\ddto\\
& &f\rrto|(0.53)\hole^(0.35){\phi}\dlto_{\pi} & & g\dlto_{\mu}\\
& h\rrto^{\psi} & & k
\enddiagram
\end{equation}
The top right square
$(\phi^{\ast}\phi_{\ast}\pi^{\ast}\sigma,\psi^{\ast}\psi_{\ast}\sigma,
\phi_{\ast}\pi^{\ast}\sigma,\psi_{\ast}\sigma)$ is a push-out by the
hypercube Lemma~\ref{lemma hypercube} since the square $(f,h,g,k)$ is a push-out.
Thus to prove the proposition it is enough to show that the
top left square
$(\pi^{\ast}\sigma,\sigma,\phi^{\ast}\phi_{\ast}\pi^{\ast}\sigma,\psi^{\ast}\psi_{\ast}\sigma)$
is also a push-out.

Assume that $\phi$ is a monomorphism. It follows that so is
$\psi$. We claim that in this case $\sigma\ra \psi^{\ast}\psi_{\ast}\sigma$ and
$\pi^{\ast}\sigma\ra \phi^{\ast}\phi_{\ast}\pi^{\ast}\sigma$ are
isomorphisms. That will be proven once we show  that  the following is a
pull-back square:
\[\xymatrix{
\sigma\rto\dto &\psi_{\ast}\sigma\dto\\
h\rto^{\psi} & k }\] 
Since the vertical morphisms of this square
are pull-backs, according to Lemma~\ref{lemma cube}, it suffices
to check that on the range level we have a pull-back of spaces:
\[\xymatrix{
R\sigma\dto\rto^(0.45){\text{id}} &R(\psi_{\ast}\sigma)\dto\\
Rh\rto^{R\psi} & Rk }\] 
As $R\psi$ is a monomorphism, this is the case.

Assume now that $\pi$ is a pull-back and a monomorphism. We use
Lemma~\ref{lemma cube}. The assumption on $\pi$ implies that the
morphisms $\pi^{\ast}\sigma\ra \sigma$,
$\phi^{\ast}\phi_{\ast}\pi^{\ast}\sigma\ra
\psi^{\ast}\psi_{\ast}\sigma$, and $\mu^{\ast}\psi_{\ast}\sigma\ra
\psi_{\ast}\sigma$ are also monomorphisms and pull-backs. Thus all
the morphisms in the top left square of the diagram~(\ref{biggd})
are  pull-backs. To see that this square is a push-out we need to
prove that it is so on the range level. Let us look at the
ranges of the top layer of the diagram~(\ref{biggd}). It is a
commutative diagram of  spaces of the form:
\[\xymatrix{
 & A\dlto\rrtou|{\text{id}} \rto& C\rto\dlto & A\dlto\\
 B\rto\rrtod|{\text{id}} & D\rto & B
}\]
where the right square $(C,D,A,B)$ is a push-out and the
diagonal maps are mono\-morphisms. We can now use the ``two out of
three" Lemma~\ref{lemma 2outof3} to conclude that the left
square $(A,B,C,D)$  is also a push-out.
\end{proof}

\section{Homotopy theory of maps}

\begin{point}\label{par prop1}
The category \mbox{$\Arrows$} can be given a model
category structure where:
\begin{itemize}
\item  a morphism $\phi$ in  \mbox{$\Arrows$} is a weak equivalence
 (cofibration)
if $D\phi$ and $R\phi$ are weak equivalences (cofibrations) in
$Spaces$;
\item a morphism $\phi:f\ra g $ in  \mbox{$\Arrows$} is a fibration
 if both
$R\phi$ and the map $Df\ra\text{lim}(Rf\xra{R\phi}
Rg\xla{g}Dg\big)$, induced by $D\phi$ and $f$, are fibrations in
$Spaces$.
\end{itemize}
\end{point}

\begin{point}\label{par prop}
The category \mbox{$\Arrows$} also supports a canonical simplicial structure:
\begin{itemize}
\item for a space $K$ and a map $f$, $f\otimes K:=f\times id_{K}$;
\item  the mapping space $\text{map}(f,g)$ is given by:
\[\text{lim}\big(\text{map}(Df,Dg)\xra{g_{\ast}}\text{map}(Df,Rg\big)\xla{f^{\ast}}\text{map}(Rf,Rg)\big).\]
 \end{itemize}
 The description of the mapping spaces is straightforward from the adjunction property \cite[II.1.3]{MR0223432}.
 This simplicial structure is compatible with the model category structure defined above (the axiom SM7 is fulfilled),
 so that the category \mbox{$\Arrows$} actually is a simplicial model category.
 \end{point}

\begin{point}
Let $F:I\ra \mbox{$\Arrows$}$ be a functor. By the universal
properties,   the morphisms $\text{hocolim}_{I}F$ and
$\text{holim}_{I}F$ are respectively naturally isomorphic to the objects in
$\text{Ho}(\mbox{$\Arrows$})$ represented by
$\text{hocolim}_{I}\pi F:\text{hocolim}_{I}DF\ra
\text{hocolim}_{I}RF$ and $\text{holim}_{I}\pi F:\text{holim}_{I}DF\ra \text{holim}_{I}RF$.
\end{point}

\section{Fiberwise decomposition}

\begin{point}
Let  $f:X\ra Y$ be a map and $\sigma:\Delta[n]\ra Y$ be a simplex.
Define $df(\sigma)\ra \Delta[n]$ to be the map that fits into the
following pull-back square in $Spaces$:
\[\xymatrix{
df(\sigma)\rto \dto& X\dto^{f}\\
\Delta[n]\rto^{\sigma} & Y }\] These maps fit into a functor
$df:Y\ra\mbox{$\Arrows$}$, indexed by the simplex category of $Y$
(see~\cite[Definition~6.1]{MR2002k:55026}). The morphisms
$\{df(\sigma)\ra f\}_{\sigma\in Y}$, given by the above
commutative squares, satisfy  the universal property of the
colimit and so $\text{colim}_{Y}df= f$.
\end{point}

Functors of the form $df$ are not arbitrary, they satisfy a certain homotopy invariance property.

\begin{definition}\label{def pseudo}
 A functor $F: I \ra \Arrows$  (indexed by a small category) is called \emph{pseudo-cofibrant}
if the morphism $\text{\rm hocolim}_{I} F \ra \text{\rm colim}_{I} F$ is an isomorphism in
$\text{Ho}(\mbox{$\Arrows$})$.
\end{definition}

Although the next proposition already appeared in Dror Farjoun's book
\cite[p.183]{dror:book}, we offer a proof illustrating the ideas
of the present paper.

\begin{proposition}\label{prop projdiag}
For any $f:X\ra Y$, the functor $df:Y\ra \mbox{$\Arrows$}$ is
pseudo-cofibrant.
\end{proposition}

\begin{proof}
Let $\mathcal{S}$ be the class of spaces $Y$ for which the
proposition is true. To prove the proposition it is enough to show
that $\mathcal{S}$ satisfies the following properties (which then
imply that $\mathcal{S}$ consists of all spaces):
\begin{enumerate}
\item $\Delta[n]\in\mathcal{S}$;
\item $\coprod Y_{i}\in \mathcal{S}$ if $Y_{i}\in\mathcal{S}$;
\item
$\text{colim}(Y_{0}\la Y_{1}\mono Y_{2})\in \mathcal{S}$ if
$Y_{i}\in \mathcal{S}$ and $Y_{1}\mono Y_{2}$ is a cofibration.
\end{enumerate}

Since $id:\Delta[n]\ra \Delta[n]$ is the terminal object of the simplex category
of $\Delta[n]$, statement (1) follows from cofinality properties of homotopy colimits (\cite[Theorem XI.9.2]{MR51:1825}
and~\cite[Theorem 30.5]{MR2002k:55026}).

Statement (2) is clear. To prove (3), set
$Y:=\text{colim}(Y_{0}\la Y_{1}\mono Y_{2})$. Let $f:X\ra Y$ be a
map. Define $f_{i}:X_{i}\ra Y_{i}$ to be the map that fits into
the following pull-back square:
\[\xymatrix{
X_{i}\rto\dto_{f_{i}} & X\dto^{f}\\
Y_{i}\rto & Y
}\]
These maps form a natural transformation between the following push-outs, where
the indicated maps are cofibrations:
\[\xymatrix{\
X\dto_{f}\ar @{}[r]|(0.46)= &
*{\text{\rm colim}\hspace{2mm}\big(\hspace{-20pt}} & X_{0}\dto_{f_{0}}
&X_{1}\lto\rmono\dto_{f_{1}} & X_{2}\dto^{f_{2}} & *{\hspace{-20pt}\big)}\\
Y \ar @{}[r]|(0.46)=& *{\text{\rm colim}\hspace{2mm}\big(\hspace{-20pt}} &
Y_{0} & Y_{1}\lto\rmono & Y_{2} & *{\hspace{-20pt}\big)} }
\]

Let $Qf\ra df\epsilon$ be a cofibrant replacement in
$\text{Fun}^{b}\big(N(Y),\mbox{$\Arrows$}\big)$ (\cite[Theorem 13.1]{MR2002k:55026})
of the composition of the diagram $df:Y\ra \mbox{$\Arrows$}$
with the forgetful functor $\epsilon:N(Y)\ra Y$,
$(\sigma_{n}\ra \cdots\ra\sigma_{0})\mapsto \sigma_{0}$
(\cite[Definition 6.6]{MR2002k:55026}). Since $N(Y_{i})\ra N(Y)$
is reduced (\cite[Example 12.10]{MR2002k:55026}
and~\cite[Proposition 5.1]{MR2002k:55026}), $Qf$ restricted to
$N(Y_{i})$ is a cofibrant replacement in
$\text{Fun}^{b}\big(N(Y_{i}),\mbox{$\Arrows$}\big)$ of the
composition of $df_{i}$ with the forgetful functor $N(Y_{i})\ra Y_{i}$.
The spaces $Y_{i}$ are assumed to be in
$\mathcal{S}$ and hence  the morphism $\text{colim}_{N(Y_{i})}Qf\ra f_{i}$
is a weak equivalence. Property (3) follows now from the basic
homotopy invariance of push-outs (\cite[Proposition 2.5.(2)]{MR2002k:55026}).
\end{proof}

\section{Homotopy pull-backs}

\begin{point}
Recall that a morphism $f\ra g$ in  \mbox{$\Arrows$} is called a
homotopy pull-back if, for some (equivalently any) weak
equivalence \mbox{$\psi:g\xra{\simeq} h$} with $h$ a fibration in
$Spaces$, the morphism $f\ra \psi_{\ast}f$ is a weak equivalence.

If $\phi:f\ra g$ is a pull-back and either $g$ or $R\phi:Rf\ra Rg$
is a fibration, then $\phi$ is a homotopy pull-back.

A homotopy pull-back $\sigma\ra f$ for which $R\sigma$ is
contractible, is called a \emph{homotopy fiber} of $f$. If
$\sigma\ra f$ and $\tau\ra f$ are homotopy fibers of
$f$ such that the images of $R\sigma$ and $R\tau$ in $Rf$ lie in the same
connected component, then $\sigma$ and $\tau$ are weakly equivalent.
\end{point}

\begin{point}\label{point bphopul}
Here is a list of some basic properties of homotopy pull-backs:
\begin{enumerate}
\item Right properness: If $\phi:f\ra g$ is a weak equivalence, then it is a homotopy pull-back.
\item Fiber characterization: A morphism $\phi:f\ra g$ is a homotopy pull-back if and
only if it induces a weak equivalence of
homotopy fibers, i.e. for any commutative square:
\[\xymatrix{
\sigma\dto\rto^{\pi} & \tau\dto\\
f\rto^{\phi} & g
}\]
if $\sigma\ra f$ and $\tau\ra g$ are respectively homotopy fibers of $f$ and $g$, then $\pi$
is a weak equivalence.
\item Two out of three: Let $\phi:f\ra g$ and $\psi:g\ra h$ be morphisms. Assume that $\psi$ is a
homotopy pull-back. Then $\phi$ is a homotopy pull-back if and
only if  $\psi\phi$ is so. Assume that $\phi$ is a homotopy
pull-back and $R\phi: Rf\ra Rg$ induces an epimorphism on
the sets of connected components. Then $\psi$ is a homotopy pull-back if and
only if $\psi\phi$ is so.
\item Disjoint union: Let $\{f_{i}\}_{i\in I}$ and
$\{g_{j}\}_{j\in J}$ be collections of maps, $h:I\ra J$ a map of
sets, and $\{\phi_{i}:f_{i}\ra g_{h(i)}\}_{i\in I}$ a collection
of homotopy pull-backs. Then the following induced morphism is
also a homotopy pull-back:
\[\coprod_{i\in I} \phi_{i}:\coprod_{i\in I}f_{i}\ra \coprod_{j\in J}g_{j}\]
\end{enumerate}
\end{point}

Note that Example~\ref{ex not23} also illustrates the failure of the full two out of three property for homotopy pull-backs.

A fibration $f$ possesses the property that for any simplex $\sigma \in Rf$ the 
morphism $df(\sigma) \ra f$ is a homotopy pull-back. Maps with that property play a crucial role in this paper.

\begin{definition}\label{pt qfibration}
A map $f:X\ra Y$ is called a \emph{quasi-fibration} if for any
morphism $\alpha:\sigma\ra \tau$ in $Y$, the morphism
$df(\alpha):df(\sigma)\ra df(\tau)$ is a weak equivalence.
\end{definition}

Even if quasi-fibrations lack the global
lifting properties enjoyed by fibrations, the local information
given by the preimages of simplices still allows to recover the
homotopy fiber.

\begin{proposition}
\label{prop weakisstrong}
A map $f$ is a quasi-fibration if and only if the morphism
$df(\sigma)\ra f$ is a homotopy pull-back for any simplex
$\sigma:\Delta[n]\ra Rf$.
\end{proposition}

\begin{proof}
If the morphisms $df(\sigma)\ra f$ are homotopy pull-backs, then
they are homotopy fibers of $f$. Thus by the homotopy invariance
of homotopy fibers, $f$ is a quasi-fibration.

Assume that $f$ is a quasi-fibration. Factor the morphism $df(\sigma)\ra f$ as
a composition $df(\sigma)\ra p\ra f$ where $Rp$ is contractible
and $ p\ra f$ is a fibration and a pull-back (hence a homotopy
pull-back). It follows that $p$ is a quasi-fibration and
$df(\sigma)\ra p$ is a pull-back. Since $Rp$ is a contractible,
according to~\cite[Lemma 27.8]{MR2002k:55026}, the morphism
$df(\sigma)\ra \text{hocolim}_{Rp}dp$ is an isomorphism in
$\text{Ho}(\mbox{$\Arrows$})$. Thus by
Proposition~\ref{prop projdiag}, $df(\sigma)\ra \text{colim}_{Rp}dp=p$ is a weak
equivalence and therefore a homotopy pull-back. We can conclude
that the composition $df(\sigma)\ra p\ra f$ is also a homotopy
pull-back.
\end{proof}

The last proposition combined with the fiber characterization of homotopy pull-backs
(property~\ref{point bphopul}.(2)) gives:

\begin{corollary}\label{cor hopullfiber}
A morphism $\phi:f\ra g$ between quasi-fibrations $f$ and $g$ is a homotopy pull-back if and only if, for any
simplex $\sigma\in Rf$, the induced morphism
$\phi(\sigma):df(\sigma)\ra dg((R\phi)\sigma)$ is a weak
equivalence.   \hfill{\qed}
\end{corollary}

\section{Distinguished collections}
In this section we prove some fundamental properties of distinguished collections
(see Definition~\ref{def dist} in the introduction). We start with a stronger form of condition (3).

\begin{proposition}\label{prop composition}
Let $\C$ be a distinguished class. If $\phi:f\ra g$ and $\psi:g\ra h$ belong to 
$\C$, then so does $\psi \phi$.
\end{proposition}

\begin{proof}
The homotopy colimit of the following diagram $F$ in $\Arrows$:
$$
f \xra{\phi} g \xla{\text{\rm id}_g} g \xra{\psi} h
$$
is homotopy equivalent to $h$. Moreover the morphism $f \ra \text{\rm hocolim} F$
can be identified with the composition $\psi \phi$. Since all morphisms in this 
diagram belong to $\C$, by condition (4)
in Definition~\ref{def dist} so does $\psi \phi$.
\end{proof}

We next characterize the collection of homotopy pull-backs.

\begin{theorem}\label{thm mainpullback}
The collection of homotopy pull-backs is the smallest
distinguished class of morphisms in \mbox{$\Arrows$}.
\end{theorem}

\begin{proof}
We first show that homotopy pull-backs form a distinguished class.
According to~\ref{point bphopul}, the requirements (1), (2), and (3) of
Definition~\ref{def dist} are satisfied. We need to prove that
homotopy pull-backs satisfy also requirement (4). We first show a
particular case:

\begin{lemma}\label{lem pushpull}
Let the following be a push-out square in \mbox{$\Arrows$}:
\[\xymatrix{
f\rto^{\phi}\dto_{\pi} &g\dto^{\mu}\\
h\rto^{\psi} & k }\] Assume that $\phi$ and $\pi$ are homotopy
pull-backs and one of them is a monomorphism. Then $\mu$ and
$\psi$ are also homotopy pull-backs.
\end{lemma}

\begin{proof}[Proof of Lemma~\ref{lem pushpull}]
By making various factorizations we may assume that both morphisms
$\phi$ and $\pi$ are monomorphisms and the maps $f$, $g$, $h$ are
fibrations.

Choose a simplex in $Rk$. Since it is in the image of either $R\mu$ or
$R\psi$, by symmetry, we can assume that it is of the form
$\Delta[n]\xra{\sigma}Rh\xra{R\psi} Rk$. It follows that
$dk(\sigma)= \psi_{\ast}dh(\sigma)$ and the following is a
push-out square (see Proposition~\ref{prop pushfiber}):
\[\xymatrix{
\pi^{\ast}dh(\sigma)\rto\dto & \phi_{\ast}\pi^{\ast}dh(\sigma)\dto\\
dh(\sigma)\rto &dk(\sigma) &
 \hspace{-14mm} = \psi_{\ast}dh(\sigma)}\]
Since $\phi$ is a homotopy pull-back between fibrations, according
to Proposition~\ref{prop projdiag} and
Corollary~\ref{cor hopullfiber}, the morphism $\pi^{\ast}dh(\sigma)\ra
\phi_{\ast}\pi^{\ast}dh(\sigma)$ is a weak equivalence.
Same is therefore true for $dh(\sigma)\ra dk(\sigma)$. We conclude that $k$ is a
quasi-fibration and $\psi$ is a homotopy pull-back
(Corollary~\ref{cor hopullfiber}).
\end{proof}

To prove in general that  homotopy pull-backs satisfy requirement
(4) of Definition~\ref{def dist}, it would be  enough to show that, for
any bounded and cofibrant $F:K\ra \mbox{$\Arrows$}$ which sends
morphisms in $K$ to homotopy pull-backs, the morphism
$F(\sigma)\ra \text{colim}_{K}F$ is a homotopy pull-back for any
simplex $\sigma\in K$.   Induction on the dimension of $K$ seems to be the right
strategy to do that. Unfortunately the notion of cofibrancy of bounded functors
is too rigid for that: cofibrant functors are not preserved by restricting along
maps of simplicial sets. To circumvent this problem we need to allow ``more general
cofibrant" diagrams. We are going to apply the idea of relative boundedness and
cofibrancy introduced in~\cite[Sections 17, 19]{MR2002k:55026} to deal with such problems.

Fix a space $L$ and denote by $\mathcal{S}_L$ the class of maps
of the form $f: K \ra L$ that satisfy
the following property: If $F~:K \ra \Arrows$ is any $f$-bounded and
$f$-cofibrant diagram which sends
morphisms in $K$ to homotopy pull-backs, then the morphism
$F(\sigma) \ra \text{colim}_K F$ is a
homotopy pull-back for any simplex $\sigma \in K$.

We will show that $\mathcal{S}_L$ satisfies the following
properties, and thus consists of all maps with range $L$:
\begin{enumerate}
\item All maps of the form $\Delta[n] \ra L$ belong to $\mathcal{S}_L$;
\item If a set of maps $K_{i} \ra L$ belongs to $\mathcal{S}_L$ then so does $\coprod K_{i} \ra L$;
\item Let the following be a commutative diagram where the indicated map is a cofibration:
\[\xymatrix{\
K\dto_{f}\ar @{}[r]|(0.46)= &
*{\text{\rm colim}\hspace{2mm}\big(\hspace{-20pt}} & K_{0}\dto_{f_{0}}
&K_{1}\lto\rmono\dto_{f_{1}} & K_{2}\dto^{f_{2}} & *{\hspace{-20pt}\big)}\\
L \ar @{}[r]|(0.46)=& *{\text{\rm colim}\hspace{2mm}\big(\hspace{-20pt}} &
L & L\lto_{\text{\rm id}}\rto^{\text{\rm id}} & L & *{\hspace{-20pt}\big)} }
\]
If each $f_i$ belongs  to $\mathcal{S }_L$ then so does $f$.
\end{enumerate}
Property (1) is clear, since $\Delta[n]$ has a terminal object.
Property (2) is easily verified as a simplex of $\coprod K_i$ is a
simplex  of one of the spaces~$K_i$. It remains to prove (3).
Let $F:K\ra \mbox{$\Arrows$}$ be an
$f$-bounded and $f$-cofibrant functor. Consider the following
commutative diagram:
\begin{equation}\label{eq pushouthpl}
\diagram
\coprod_{\sigma\in K_{1}}F(\sigma)\rrto^{b}\drto^{d}\ddto_{a}
& & \coprod_{\sigma\in K_{2}}F(\sigma)\dto^{e}\\
 &\text{colim}_{K_{1}}F\rmono\dto & \text{colim}_{K_{2}}F\dto\\
\coprod_{\sigma\in K_{0}}F(\sigma)\rto^{c} &
\text{colim}_{K_{0}}F\rmono & \text {colim}_{K}F
\enddiagram\end{equation}
According to the disjoint union property~\ref{point bphopul}.(4),
the morphisms $a$ and $b$ are homotopy pull-backs. The functor  $F$
restricted along $f_{i}:K_i\ra L$ is both  $f_i$-bounded and $f_i$-cofibrant
(see~\cite[Corollaries 17.5 (1), 19.6 (1)]{MR2002k:55026}). Thus
by the inductive hypothesis $c$, $d$, and $e$ are homotopy
pull-backs. We can now apply the two out of three
property \ref{point bphopul}(3) to see that both morphisms
$\text{colim}_{K_1} F \mono \text{colim}_{K_2}F$ and
$\text{colim}_{K_1} F \ra \text{colim}_{K_0} F$ are homotopy pull-backs too. By
Lemma~\ref{lem pushpull} we can conclude that all the morphisms in the
diagram~(\ref{eq pushouthpl}) are homotopy pull-backs.

We are left to show that homotopy pull-backs are contained in any
distinguished class. For that it is enough to show that any
pull-back $\phi: f\ra g$ with $g$ a fibration belongs to any
distinguished class.

Assume first that  $\phi$ coincides with $dg(\sigma)\ra g$, for
some $\sigma\in Rg$. Note that the functor $dg:Rg\ra \mbox{$\Arrows$}$
is pseudo-cofibrant (Proposition~\ref{prop projdiag}) and it takes
all the morphisms in $K$ to weak
equivalences. As weak equivalences belong to any distinguished
class, then so does $dg(\sigma)\ra g$.

For a general pull-back $\phi:f\ra g$, define $I$ to be the
Grothendieck construction
$I:=\text{Gr}(\Delta[0]\la Rf\xra{R\phi}Rg)$
(\cite[Section 38]{MR2002k:55026}). Define
further  $F:I\ra \mbox{$\Arrows$}$ to be the functor given by the
data (see \cite[Section 40]{MR2002k:55026}):
\begin{equation}
\label{eq groth}
\begin{array}{l}
F_{\Delta[0]}:\Delta[0]\ra \mbox{$\Arrows$} \text{ is the constant functor
with value } f;\\
F_{Rf}:  Rf\ra \mbox{$\Arrows$} \text{ is } df;\\
F_{Rg}: Rg\ra \mbox{$\Arrows$} \text{ is } dg;\\
F_{Rf}\ra F_{\Delta[0]} \text{ is given by the morphisms } df\ra\text{colim}_{Rf}df=f;\\
F_{Rf}\ra F_{Rg} \text{ is induced by } \phi.
\end{array}
\end{equation}
Again by Proposition~\ref{prop projdiag}, the functor $F$ is
pseudo-cofibrant. Moreover it takes any morphism in $I$ to either
a weak equivalence or a morphism of the form $df(\sigma)\ra f$.
Since such morphisms belong to any distinguished class we conclude
that so does $\phi:f=F(\Delta[0])\ra \text{colim}_{I}F=g$.
\end{proof}

Note that we have not used condition (3) in Definition~\ref{def dist}
while proving Theorem~\ref{thm mainpullback}.
This means that the collection of homotopy pull-backs can be
characterized as the smallest class that satisfies only the three
other requirements of the definition. The significance of the
third condition is illustrated by:

\begin{corollary}\label{cor fiberchar}
Let $\C$ be a distinguished class. Then $\phi: f \ra g$ belongs to $\C$
if and only if the morphism $\pi: \sigma \ra \tau$
does so for any commutative square of the form:
\[\xymatrix{
\sigma\dto\rto^{\pi} & \tau\dto\\
f\rto^{\phi} & g
}\]
where $\sigma\ra f$ and $\tau\ra g$ are respectively
homotopy fibers of $f$ and $g$.
\end{corollary}

\begin{proof}
Since $\sigma \ra f$ and $\tau \ra g$ are homotopy pull-backs they
belong to any distinguished class by
Theorem~\ref{thm mainpullback}. Let us assume that $\phi$ is a member of $\C$. By
Proposition~\ref{prop composition} the morphism $\sigma \ra g$
is in $\C$ and hence by condition (3) in Definition~\ref{def dist} so is $\pi$.

Let us prove now the converse. By making an appropriate
factorization we can assume that $f$ and $g$ are fibrations.
Let $F: I \ra \Arrows$ be the functor given by (\ref{eq groth}) in the proof of 
Theorem~\ref{thm mainpullback}.
This functor takes any morphism in $I$ either to a morphism
in $\C$ (by assumption) or to a homotopy pull-back,
which is also in $\C$. Therefore
$\phi:f=F(\Delta[0])\ra \text{colim}_{I}F=g$ belongs to $\C$.
\end{proof}

\section{Fiberwise localization}\label{sec fiblocal}

\begin{point}
Let $\phi$ be a morphism in \mbox{$\Arrows$}.
 Recall that a  map of spaces $f$ is called $\phi$-local if, for some (equivalently any)
weak equivalence $f\simeq g$ with $g$ a fibrant, the map of spaces
$\text{map}(\phi,g)$ is a weak equivalence.

According to~\cite{MR0478159}, ~\cite{MR95c:55010},~\cite{dror:book},
and~\cite{MR1944041}, there is a functor
$L_{\phi}:\mbox{$\Arrows$}\ra \mbox{$\Arrows$}$ and a natural
transformation $f\ra L_{\phi}f$ (called localization) such that:
\begin{itemize}
\item $L_{\phi}f$ is fibrant and $\phi$-local;
\item  the map of spaces
$\text{map}(L_{\phi}f,g)\ra \text{map}(f,g)$ is a weak equivalence,
for any fibrant and $\phi$-local $g$.
\end{itemize}

We are going to refer to $\phi$-local maps also as
$L_{\phi}$-local and to morphisms $\psi$ for which $L_{\phi}\psi$
is a weak equivalence as  $L_{\phi}$-equivalences.
\end{point}

In general it is very difficult to understand the
localization with respect to an arbitrary
morphism $\phi$. However when $\phi$ is a morphism of the
form $(u,\text{id}_{\Delta[0]})$
the next two propositions identify the localization with a very
familiar object. Through out this section we are going to fix a map of spaces $u: A\ra B$.
The  symbol $L_{u}$ will be used  to denote both the
localization in $Spaces$ with respect to $u$ and the localization
in \mbox{$\Arrows$} with respect to the morphism $(u,\text{id}_{\Delta[0]})$. 
We will   show that the functor $L_{u}$ in  \mbox{$\Arrows$}
is the fiberwise version of $L_{u}$ in $Spaces$.

We start by  characterizing the $L_u$-local maps as those for
which the homotopy fibers are $L_{u}$-local spaces.

\begin{proposition}\label{prop fibrlocal}
 A map $f$ is $L_{u}$-local in $\Arrows$ if and only if,
for any homotopy fiber $\sigma\ra f$, the space $D\sigma$ is $L_{u}$-local
in $Spaces$, i.e. for any weak equivalence $D\sigma\simeq Z$ with
$Z$ fibrant, the map of spaces $\text{\rm map}(u,Z)$ is a weak equivalence.
\end{proposition}

\begin{proof}
Choose a weak equivalence $f\simeq g$ with $g$  fibrant ($Rg$ is fibrant and
$g$ is a fibration). By definition of the simplicial structure given in ~\ref{par prop}
we have the following cube of spaces :
\[\xymatrix{
&\text{map}(B\ra\Delta[0],g)\rrepi^{b}\dlto \ddto|(0.26){\text{map}(\phi,g)}|\hole
& &Rg\dlto \ddto|{\text{id}} \\
\text{map}(B,Dg)\rrepi^(0.3){\text{map}(B,g)} \ddto|{\text{map}(u,Dg)}
& & \text{map}(B,Rg) \ddto|(0.26){\text{map}(u,Rg)} \\
&\text{map}(A\ra\Delta[0],g)\ar@{->>}[rr]|(0.61)\hole^{a} \dlto & &Rg\dlto\\
\text{map}(A,Dg)\rrepi^{\text{map}(A,g)} & & \text{map}(A,Rg)
}\]
where the top and bottom faces are pull-back squares and the labeled arrows are 
fibrations. Let $x\in Rg$ be a vertex.
The fibers of $a$ and $b$ over the vertex $x$
can be identified respectively with the mapping spaces
$\text{map}(A,g^{-1}(x))$ and $\text{map}(B,g^{-1}(x))$.
Thus $\text{map}(\phi,g)$
is a weak equivalence if and only if, for any vertex $x\in Rg$,
the map $\text{map}(u,g^{-1}(x))$
is a weak equivalence.
\end{proof}

In general local objects are not closed under homotopy colimits. However in the 
case of $L_{u}$ we have:

\begin{corollary}\label{col localcolim}
Let $F:I\ra \mbox{$\Arrows$}$ be a
pseudo-cofibrant functor. Assume that, for any $i$, the map $F(i)$
is $L_{u}$-local and, for any morphism $\alpha$ in $I$, the morphism
$F(\alpha)$ is a homotopy pull-back. The map $\text{\rm colim}_{I}F$
is then $L_{u}$-local.
\end{corollary}

\begin{proof}
Since for any $i\in I$, $F(i)\ra \text{\rm colim}_{I}F$ is a
homotopy pull-back (see Theorem~\ref{thm mainpullback}), any
homotopy fiber of $\text{\rm colim}_{I }F$ is a homotopy
fiber of some $F(i)$. As these  are $L_{u}$-local spaces,
the map $\text{\rm colim}_{I}F$ is $L_{u}$-local in
\mbox{$\Arrows$}.
\end{proof}

Next we  describe $L_u$-equivalences:

\begin{proposition}\label{prop h equivalence}
 A morphism $\psi$  in
\mbox{$\Arrows$} is an  $L_{u}$-equivalence if and only if:
\vskip -2mm
\vbox{\begin{itemize}
\item
$R\psi$ is a weak equivalence;
\item for any commutative square:
\[\xymatrix{
\sigma\rto^{\pi}\dto & \tau\dto\\
f\rto^{\psi} & g
}\]
if $\sigma\ra f$ and $\tau\ra g$ are homotopy fibers, then
$D\pi$ is an $L_{u}$-equivalence in $Spaces$.
\end{itemize}}
\end{proposition}

\begin{proof}
Assume first that $\psi:f\ra g$ is an
$L_{u}$-equivalence. Note that for any fibrant space $X$, the map
$\text{id}_{X}$ is $L_{u}$-local by Proposition~\ref{prop fibrlocal}.
It follows that
$\text{map}(\psi, \text{id}_{X})= \text{map}(R\psi, X)$
is a weak equivalence of spaces for all fibrant $X$. Hence $R\psi$
is a weak equivalence.

If $Z$ is a fibrant and $L_{u}$-local space, then the map $Z\ra \Delta[0]$ is  $L_{u}$-local
 in \mbox{$\Arrows$}. Thus the map of spaces:
 \[\xymatrix{
 \text{map}(Dg,Z)\dto_{\text{map}(D\psi ,Z)}\ar@{}[r]|(0.44)=
 & \text{map}(g,Z \ra\Delta[0])\dto^{\text{map}(\psi,Z\ra \Delta[0])}\\
 \text{map}(Df,Z) \ar@{}[r]|(0.44)= & \text{map}(f,Z\ra\Delta[0])
 }\]
is a weak equivalence. This together with the fact that $R\psi$
is a weak equivalence implies that, for any commutative square:
 \[\xymatrix{
 \sigma\dto\rto^{\pi} & \tau\dto\\
 f\rto^{\psi} &g
 }\]
where $\sigma\ra f$ and $\tau\ra g$ are homotopy fibers, then
 $\text{map}(D\pi,Z)$ is a weak equivalence of spaces. As this holds for any
 $L_{u}$-local space $Z$, the map $D\pi$ is an $L_{u}$-equivalence.

To prove the other implication consider the following commutative square:
\[\xymatrix{
f\rto^{\psi}\dto & g\dto\\
L_{u}f\rto ^{L_{u}\psi} & L_{u}g }\] 
The vertical morphisms are
$L_{u}$-equivalences and we already know they
induce weak equivalences on ranges and $L_{u}$-equivalences on
homotopy fibers. Thus if $\psi$ satisfies the two properties of
the proposition, then so does $L_{u}\psi$. Since $L_{u}\psi$ is a
morphism between $L_{u}$-local maps, i.e. maps whose homotopy
fibers are $L_{u}$-local spaces (see
Proposition~\ref{prop fibrlocal}), the morphism
$L_{u}\psi$ is a weak equivalence. We
can conclude that $\psi$ is an $L_{u}$-equivalence.
\end{proof}

\section{$L_{u}$-homotopy pull-backs}\label{sec luhpullback}

In Theorem~\ref{thm mainpullback} we saw that any
distinguished class containing all weak equivalences must
also contain all homotopy pull-backs. The analogous statement
for $L_\phi$-equivalences should involve
$L_\phi$-homotopy pull-backs.

\begin{definition}\label{def h-homotopypb}
{\rm A morphism $f \ra g$ in $\Arrows$ is called an
\emph{$L_\phi$-homotopy pull-back} if, for
some (equivalently any) weak equivalence
$\psi: g \stackrel{\simeq}{\ra} h$ with $h$ a fibration
in $Spaces$, the morphism $f \ra \psi_* f$ is an
$L_\phi$-equivalence.
}
\end{definition}

There is a more amenable description of $L_\phi$-homotopy pull-backs when $\phi$
 is of the form $(u,\text{id}_{\Delta[0]})$.

\begin{proposition}\label{prop amenable}
If $u$ is a map of spaces, then $\pi : f \ra g$   is an $L_u$-homotopy pull-back
 if and only if $L_u \pi$ is a homotopy pull-back.
\end{proposition}

\begin{proof}
From the fiber characterization property~\ref{point bphopul}.(2) we see that the
morphism $L_u\pi~: L_u f \ra L_u g$ is a homotopy pull-back if and only if  it
induces a weak equivalence on homotopy fibers. According to
Proposition \ref{prop h equivalence} this
is the case if and only if the morphism $\pi$ induces $L_u$-equivalences on
homotopy fibers.
For any weak equivalence $\psi: g \stackrel{\simeq}{\ra} h$ with $h$ a fibration
in $Spaces$,
the morphism $\psi_*f \ra  h$ always induces an equivalence on homotopy fibers. 
By construction, $f$ and $\psi_* f$ have the same range. Therefore, the
morphism $f \ra \psi_* f$ is an $L_u$-equivalence if and only if
$L_u\pi$ is a homotopy pull-back.
\end{proof}

Using this characterization, we can now prove our main result:

\begin{theorem}\label{thm luhpullback}
Let $u$ be a map of spaces. The collection of $L_u$-homotopy pull-backs is the smallest distinguished
class  containing all $L_{u}$-equivalences.
\end{theorem}

\begin{proof}
Consider a commutative square:\[\xymatrix{
f\rto^{\pi}\dto & g\dto\\
L_{u}f\rto^{L_{u}\pi} & L_{u}g
}\]
If $\pi$ is an $L_u$-homotopy pull-back, then $L_{u}\pi$ is
a homotopy pull-back and hence by
Theorem~\ref{thm mainpullback} it belongs to any
distinguished class. By condition (3)
of Definition~\ref{def dist} it follows that $\pi$ belongs to any distinguished 
class that contains $L_{u}$-equivalences.

To prove the theorem it remains to show that the collection in the statement is a
distinguished class. It is clear that requirements (1), (2), and (3) of
Definition~\ref{def dist} are satisfied. Let $F:I\ra
\mbox{$\Arrows$}$ be a pseudo-cofibrant functor which takes
morphisms in $I$ to $L_u$-homotopy pull-backs. We need to show that, for
any $i\in I$, $L_{u}F(i)\ra L_{u}(\text{colim}_{I}F)$ is a
homotopy pull-back.
The functor $F$ fits into the following commutative square:
\[\xymatrix{
H\rto\dto_{\simeq} & G\dto^{\simeq}\\
F\rto & L_{u}F }\] where $H$ and $G$ are pseudo-cofibrant  and the
indicated natural transformations are weak equivalences
(see \cite[Remark 16.3]{MR2002k:55026}).  Choose an object $i\in I$ and
consider the following commutative diagram, where the indicated arrows are weak equivalences:
\[\xymatrix{
L_{u}F(i)\dto_{e} & L_{u}H(i)\lto_{\simeq}\dto\rto^{\simeq}
& L_{u}G(i)\dto & G(i)\lto_{b}\dto^{a}\\
L_{u}(\text{colim}_{I}F) &
L_{u}(\text{colim}_{I}H)\lto_{\simeq}\rto^{d} & L_{u}
(\text{colim}_{I}G) &\text{colim}_{I}G\lto_(0.4){c} }\]
As $G(\alpha)$ is a homotopy pull-back for any morphism $\alpha\in I$,
by Theorem~\ref{thm mainpullback}, the morphism $a$ is a
homotopy pull-back. Furthermore the values of $G$ are
$L_{u}$-local so $b$ is a weak equivalence.  It follows from
Corollary~\ref{col localcolim} that $\text{colim}_{I}G$ is also
$L_{u}$-local and
consequently the morphism $c$ is a weak equivalence. Since
$L_{u}$-equivalences are preserved by homotopy colimits, the
morphism $d$ is a weak equivalence. We can therefore conclude that
$e$ is a homotopy pull-back.
\end{proof}

To illustrate this theorem we offer an application with a classical flavor.
If for a map of spaces all the preimages of
simplices have the same homotopy type up to $L_u$-localization,
then the $L_u$-localization of the homotopy fiber shares the
same homotopy type as well (the proof is identical
as that of Proposition~\ref{prop weakisstrong}):

\begin{corollary}
\label{cor weakisstrong}
Let $f: X \ra Y$ be a map over a connected space $Y$. Let  $F$ be its homotopy fiber.
Then the morphisms $df(\sigma)\ra f$ are $L_u$-homotopy pull-backs for all simplices
$\sigma:\Delta[n]\ra Y$ if and only if each of the induced map $Ddf(\sigma)
\ra F$ is an $L_u$-equivalence. \hfill{\qed}
\end{corollary}

The particular case when the preimages of simplices are acyclic with respect to 
some generalized homology
theory is due to E. Dror Farjoun \cite[Corollary~9.B.3.2]{dror:book}.
\begin{corollary}
\label{cor acyclicfibers}
Let $E$ be a generalized homology theory and $f: X \ra Y$ be a map over a connected space $Y$
such that $df(\sigma)$ is $E$-acyclic for any simplex
$\sigma:\Delta[n]\ra Y$. The homotopy fiber of $f$ is then also $E$-acyclic. \hfill{\qed}
\end{corollary}

\section{Concluding remark}
According to Corollary~\ref{cor fiberchar} a distinguished collection $\C$ is determined by
the class $\text{T}(\C)$ of maps $f: A \ra B$ of spaces for which the morphism
$(f,\text{id}_{\Delta[0]})$ in $\Arrows$ belongs to $\C$. Members of  $\text{T}(\C)$ are called
\emph{transition functions} for $\C$. For example the transition functions for
homotopy pull-backs are weak equivalences and the transition functions for
$\text{H}\mathbf Z$-homotopy pull-backs are $\text{H}\mathbf Z$-equivalences.
More generally, if $u$ is a map of spaces, then  the transition functions for $L_{u}$-homotopy pull-backs
are $L_{u}$-equivalences.

In these three examples, the class $\text{T}(\C)$
satisfies  the following properties:
\begin{enumerate}
\item[(A)] $\text{T}(\C)$ contains weak equivalences.
\item[(B)] Let $f:A\ra B$ and $g:B\ra C$ be maps. If two out of $f$, $g$ and $gf$ are in $\text{T}(\C)$,
then so is the third.
\item[(C)] Let $F:I\ra \Arrows$ be a functor. If, for any $i\in I$, $F(i)$ belongs to $\text{T}(\C)$, then
so does $\text{hocolim}_{I}F$.
\end{enumerate}

On the other hand we could start with a class of maps $\D$ and define $\D$-homotopy pull-backs to
be the collection $\C(\D)$ of morphisms $\phi:f\ra g$ in $\Arrows$
for which the maps induced on homotopy fibers of $f$ and $g$  belong to $\D$.

While writing this paper we have not found answers to the following questions:

\begin{question}
\label{question1}
Is it true that, for any distinguished collection $\C$, the class of transition 
functions $\text{T}(\C)$
satisfies the above three conditions (A), (B), and (C)?
\end{question}

Note that if  $\text{T}(\C)$  satisfies the full ``two out of three" condition (B), then
$\C$ has the following ``extended  two out of three" condition:
Let  $\phi:f\ra g$ and $\psi:g\ra h$ be morphisms. Assume that $R\phi$ is an
epimorphism on the set of connected components. Then $\psi$ belongs to $\C$ if and
only if $\psi\phi$ does. This should be compared with property (3) in
Section~\ref{point bphopul}. We do not know if
any distinguished collection satisfies such an ``extended  two out of three" condition.
Maybe this extra requirement ought to be added to the
definition of a distinguished collection. In this paper
though we tried to avoid making any general connectivity assumption.
\medskip

\begin{question}
\label{question2}
What are the necessary and sufficient requirements on a class of maps $\D$, so that
$\C(\D)$ is a distinguished collection?
\end{question}

Requirements  (A), (B), and (C) above are possible candidates.

\bibliographystyle{plain}\label{biblography}
\bibliography{bibho}

\begin{thebibliography}{10}

\bibitem{MR80d:55001}
J.~F. Adams.
\newblock {\em Infinite loop spaces}, volume~90.
\newblock Princeton University Press, 1978.
\newblock Annals of Mathematics Studies.

\bibitem{MR47:3489}
M.~Baratt and S.~Priddy.
\newblock On the homology of non-connected monoids and theirs associated
  groups.
\newblock {\em Comment. Math. Helv.}, 47:1--14, 1972.

\bibitem{MR0478159}
A.~K. Bousfield.
\newblock Constructions of factorization systems in categories.
\newblock {\em J. Pure Appl. Algebra}, 9(2):207--220, 1976/77.

\bibitem{MR95c:55010}
A.~K. Bousfield.
\newblock Localization and periodicity in unstable homotopy theory.
\newblock {\em J. Amer. Math. Soc.}, 7(4):831--873, 1994.

\bibitem{MR51:1825}
A.~K. Bousfield and D.~M. Kan.
\newblock {\em Homotopy limits, completions and localizations}.
\newblock Springer-Verlag, Berlin, 1972.
\newblock Lecture Notes in Mathematics, Vol. 304.

\bibitem{MR2002k:55026}
W. Chach{\'o}lski and J. Scherer.
\newblock Homotopy theory of diagrams.
\newblock {\em Mem. Amer. Math. Soc.}, 155(736):x+90, 2002.

\bibitem{dror:book}
E.~Dror~Farjoun.
\newblock {\em Cellular spaces, null spaces and homotopy localization}, volume
  1622 of {\em Lecture Notes in Mathematics}.
\newblock Springer-Verlag, Berlin, 1996.

\bibitem{MR93i:55015}
T.~G. Goodwillie.
\newblock Calculus. {II}. {A}nalytic functors.
\newblock {\em $K$-Theory}, 5(4):295--332, 1991/92.

\bibitem{MR1944041}
P.~S. Hirschhorn.
\newblock {\em Model categories and their localizations}, volume~99 of {\em
  Mathematical Surveys and Monographs}.
\newblock American Mathematical Society, Providence, RI, 2003.

\bibitem{MR90j:55029}
J.~F. Jardine.
\newblock The homotopical foundations of algebraic ${K}$-theory.
\newblock In {\em Algebraic $K$-theory and algebraic number theory (Honolulu,
  HI, 1987)}, pages 57--82. Amer. Math. Soc., Providence, RI, 1989.

\bibitem{MR595633}
M. Mather and M. Walker.
\newblock Commuting homotopy limits and colimits.
\newblock {\em Math. Z.}, 175(1):77--80, 1980.

\bibitem{MR53:6547}
D.~McDuff and G.~Segal.
\newblock Homology fibrations and the ``group-completion'' theorem.
\newblock {\em Invent. Math.}, 31(3):279--284, 1975/76.

\bibitem{MR2061572}
W. Pitsch and J. Scherer.
\newblock Homology fibrations and ``group-completion'' revisited.
\newblock {\em Homology Homotopy Appl.}, 6(1):153--166 (electronic), 2004.

\bibitem{MR51:1808}
V. Puppe.
\newblock A remark on ``homotopy fibrations''.
\newblock {\em Manuscripta Math.}, 12:113--120, 1974.

\bibitem{MR0223432}
D.~G. Quillen.
\newblock {\em Homotopical algebra}.
\newblock Lecture Notes in Mathematics, No. 43. Springer-Verlag, Berlin, 1967.

\bibitem{MR99k:57036}
U.~Tillmann.
\newblock On the homotopy of the stable mapping class group.
\newblock {\em Invent. Math.}, 130(2):257--275, 1997.

\end{thebibliography}

\end{document}